\documentclass[11pt]{article}
\usepackage{amsmath,amsthm,latexsym,amssymb}
\usepackage{amsfonts}
\usepackage{amsbsy}
\usepackage{hyperref}

\newtheorem{theo}{\bf Theorem}[section]

\newtheorem{coro}[theo]{\bf Corollary}
\newtheorem{rem}[theo]{\bf Remark}

\def\be{\begin{eqnarray}}
\def\ee{\end{eqnarray}}

\def\RR{\mathbb{R}}
\def\NN{\mathbb{N}}
\def \EE{\mathbb{E}}
\def \PP{\mathbb{P}}

\def \1{\mathbf{1}}

\renewcommand{\Box}{\hfill\mbox{\fbox{\rule{0mm}{1.5mm}}}}

\let\la=\lambda

\begin{document}

\title{\textbf{Measurability of optimal transportation and strong coupling of martingale measures}}
\author{Joaquin Fontbona\thanks{DIM-CMM, UMI(2807) UCHILE-CNRS, Universidad de Chile, Casilla 170-3, Correo 3,
Santiago-Chile, e-mail:fontbona@dim.uchile.cl. Supported by
Fondecyt Proyect 1070743, ECOS-Conicyt C05E02, Millennium Nucleus
Information and Randomness ICM P04-069-F and FONDAP Applied
Mathematics}, H\'el\`ene Gu\'erin
\thanks{IRMAR, Universit\'e Rennes 1, Campus de Beaulieu, 35042 Rennes-France,
e-mail:helene.guerin@univ-rennes1.fr. Supported by  ECOS-Conicyt
C05E02}, Sylvie M\'el\'eard\thanks{CMAP, Ecole Polytechnique,
CNRS, route de Saclay, 91128 Palaiseau Cedex-France e-mail:
sylvie.meleard@polytechnique.edu.  Supported by ECOS-Conicyt
C05E02 and Millennium Nucleus Information and Randomness ICM
P04-069-F}}

\date{}

\pagenumbering{arabic} \maketitle
 \pagestyle{plain} \thispagestyle{empty}

\begin{abstract}
We consider  the optimal mass transportation problem in $\RR^d$
with  measurably parameterized marginals, for general cost
functions and under  conditions ensuring the existence of a unique
optimal transport map. We prove a joint measurability result for
this map, with respect to the space variable and to the parameter.
The proof needs to establish the measurability of some  set-valued
mappings, related to the support of the optimal transference
plans, which we  use to perform a suitable discrete approximation
procedure. A motivation is the construction of a strong coupling
between orthogonal martingale measures. By this we mean that,
given a martingale measure, we construct in the same probability
space  a second one with specified covariance measure. This is
done by pushing forward one martingale measure through a
predictable version of the optimal transport map between  the
covariance measures. This coupling  allows us to obtain
quantitative estimates in terms of the Wasserstein distance
between those covariance measures.


\end{abstract}

\section{Introduction}
We consider  the optimal mass transportation problem in $\RR^d$
with  measurably parameterized marginals, for general cost
functions and under  conditions ensuring the existence of a unique
optimal transport map. The aim of this note is to prove a joint
measurability result for this map, with respect to the space
variable and to the parameter. One of our motivations, developed
at the end, is the construction of a strong coupling between
martingale measures. That is, given a martingale measure, we shall
construct in the same probability space  a second one with
specified covariance measure process. This will be done by pushing
forward the given martingale measure through the optimal transport
map between the covariance measures. To make this construction
rigorous, we need the existence of a {\it predictable} version of
this transport map, which will be a consequence of our main
result.

\medskip

We denote the space of Borel probability measures in $\RR^d$ by
$\mathcal{ P}(\RR^d)$, and by $\mathcal{ P}_p(\RR^d)$ the subspace
of probability measures having finite $p-$order moment.

Given $\pi \in \mathcal{ P}(\RR^{2d})$,  we write
$$\pi<^{\mu}_{\nu}$$
if $\mu, \nu \in \mathcal{ P}(\RR^{d})$ are respectively its first
and second marginals. Such $\pi$ is refereed to as a
``transference plan'' between $\mu$ and $\nu$.

Let $c:\RR^d\to \RR_+$ be a continuous function. The mapping
$$\pi \to I(\pi):= \int_{\RR^{2d}} c(x,y)
\pi(dx,dy)$$ is then lower semi continuous.

The Monge-Kantorovich or optimal mass transportation problem with
cost $c$ and marginals $\mu$, $\nu$ consists in finding
$$\inf_{\pi<^{\mu}_{\nu}}I(\pi). $$
It is well known that the infimum is attained as soon as it is
finite,  see \cite{Villani:03}, Ch.1. In this case, we denote by
$\Pi_c^*(\mu,\nu)$ the subset of ${\cal P}(\RR^{2d})$ of
minimizers. If otherwise,  $I(\pi)=+\infty$ for all
$\pi<^{\mu}_{\nu}$, then by convention we set
$\Pi_c^*(\mu,\nu)=\emptyset$.

\medskip

We shall say that {\it Assumption $H(\mu,\nu,c)$} holds if
\begin{itemize}
\item[a)] $\mu$ {\it does not give mass to sets
 with Hausdorff dimension smaller than or equal to $d-1$.}
\item[b)] {\it there exists a unique optimal transference  plan $\pi
\in \Pi_c^*(\mu,\nu)$, and it has the form
$$\pi(dx,dy)=\mu(dx)\otimes
\delta_{T(x)}(dy)$$ for a  $\mu(dx)-a.s.$ unique mapping
$T:\RR^d\to \RR^d$.}
\end{itemize}
 Such $T$ is called an {\it optimal transport map between $\mu$ and
$\nu$ for the cost function $c$.}

 Hypothesis a) in $H(\mu,\nu,c)$ is optimal both
for existence and uniqueness of an optimal transport map, see
Remark 9.5 in \cite{Villani2}. We recall that if
$\Pi_c^*(\mu,\nu)\not=\emptyset$,
 a) implies b) in the following situations (see
Gangbo and McCann \cite{GM}):
\begin{itemize}
\item[i)] $c(x,y)=\tilde{c}(|x-y|)$  with
$\tilde{c}:\RR_+\to \RR_+$ strictly convex,
  superlinear and differentiable with locally
Lipschitz gradient.
\item[ii)] $c(x,y)=\tilde{c}(|x-y|)$ with
$\tilde{c}$ strictly concave, and $\mu$ and $\nu$ are mutually
singular.
\end{itemize}
Condition b) also holds if
\begin{itemize}
\item[iii)] $c(x,y)=\tilde{c}(|x-y|)$ with
$\tilde{c}$ strictly convex and
  superlinear, and  moreover $\mu$ is absolutely continuous with respect to
Lebesgue measure.
\end{itemize}
When $\mu,\nu\in {\cal P}_p(\RR^d)$, fundamental examples  are the
cost  function  $c(x,y)=|x-y|^p$ with $p\geq 2$ for case i), $p>1$
for case iii), and $p\in (0,1)$ for case ii).

\medskip

Our main result is

\begin{theo}\label{theo:main}
Let $(E,\Sigma,m)$ be a $\sigma-$finite measurable space and
consider a measurable function $ \lambda \in E \mapsto
(\mu_\la,\nu_\la)\in \mathcal{P}(\RR^d)$ such that  for $m-$almost
every $\lambda$, $H(\mu_{\lambda},\nu_{\lambda},c)$ holds, with
optimal transport map $T_\la:\RR^d \to \RR^d$. Then, there exists
a function $(\la,x)\mapsto T(\la,x)$ which is measurable with
respect to $\Sigma\otimes {\cal B}(\RR^d)$ and such that
$m(d\la)-$almost everywhere,
$$T(\la,x)=T_\la(x)\quad \mu_\la(dx)\mbox{-almost surely.} $$
In particular, $T_\la(x)$ is measurable with respect to the
completion of $\Sigma\otimes {\cal B}(\RR^d)$ with respect to
$m(d\la )\mu_{\la}(dx)$.
\end{theo}

Theorem \ref{theo:main} generalizes Theorem 1.2 in \cite{FGM},
where we constructed a predictable version of a quadratic
transport map, between a time-varying law and empirical samples of
it.

 To our knowledge, other measurability results on the mass
transportation problem  require a topological structure on the
space of parameters, or concern transference plans but not
transport maps (see e.g. \cite{Ru}, or Corollaries 5.22 and 5.23
in \cite{Villani2}).

The proof of Theorem \ref{theo:main} is developed in the following
section. We firstly establish  a type of measurable  dependence of
the support of the optimizers on $\lambda$. From this result, we
can define
 measurable partitions of $E\times \RR^d$ induced by a dyadic
 partition of $\RR^d$, and construct bi-measurable discrete approximations of
 $T(\la,x)$. This  approximation procedure
was not needed in the simpler case studied in \cite{FGM}, where
 one of the marginals was an empirical measure (thus with finite
support).

\section{Proof of Theorem \ref{theo:main}}

Let us first state an intermediary result concerning measurability
properties of minimizers in the general framework. Its formulation
and proof require some notions of set-valued analysis, see e.g.
Appendix A of \cite{Rocka:98}.

\begin{theo}\label{theo:suppmeasurable}
 The function
assigning to $(\mu,\nu)$ the set of $\RR^{2d}$
\begin{equation}\label{eq:psi}
\Psi(\mu,\nu) := Adh\left(\bigcup_{ \pi \in \Pi_c^*(\mu,\nu)} supp
(\pi)\right),
\end{equation}
is measurable in the sense of set-valued mappings. That is, for
any open set $\theta$ in $\RR^{2d}$,  its inverse image
$\Psi^{-1}(\theta)=\{(\mu,\nu)\in ({\cal P}(\RR^d))^2:
\Psi(\mu,\nu)\cap \theta \not= \emptyset\}$ is a Borel set in
$({\cal P}(\RR^d))^2$.
\end{theo}

\begin{rem}
In the case of a set-valued mapping taking closed-set values,
measurability is equivalent to the fact that inverse images of
closed sets are measurable (see \cite{Rocka:98}).
\end{rem}

{\bf Proof.} The idea of the proof is similar to the one of
Theorem 1.3 in \cite{FGM}, where we considered the quadratic cost
and the measurable structure induced by the Wasserstein topology.
In the present case, the spaces ${\cal P}(\RR^d)$ and ${\cal
P}(\RR^{2d})$ are endowed with the usual weak topology.

 We observe that $\Psi$ writes as the adherence of a set-valued
 composition,
$$\Psi(\mu,\nu)=Adh\left(U\circ S(\mu,\nu)\right):=Adh\left(\bigcup_{\pi\in S(\mu,\nu)} U(\pi)\right),$$ where $S$ and $U$ are the
set-valued mappings respectively defined by $$S(\mu,\nu):=
\Pi_c^*(\mu,\nu) \ \mbox{ and } \ U(\pi):=supp(\pi).$$
Measurability of $\Psi$ is equivalent to $U\circ S$ being
measurable. The latter will be true as soon as $S$ is measurable
and $U^{-1}(\theta)$ is open for every open set $\theta$ (see
\cite{Rocka:98}).

The stability  theorem for optimal transference plans of
Schachermayer and Teichman (Theorem 3 in \cite{ST}) exactly states
that  inverse images through $S$ of closed sets in ${\cal
P}(\RR^{2d})$ are closed sets in $({\cal P}(\RR^d))^2$. This,
together with the fact that the mapping $S$ takes closed-set
values (by lower semi continuity of $I(\pi)$) imply that $S$ is a
measurable multi-application.

On the other hand, the inverse image by $U$ of an open set
$\theta$ of $\RR^{2d}$  is
$$U^{-1}(\theta)= \{\pi \in {\cal P}(\RR^{2d}): supp(\pi) \cap \theta \not =
\emptyset\} = \{\pi \in {\cal P}(\RR^{2d}): \pi(\theta)>0 \}. $$
It then follows by the Portmanteau Theorem that $U^{-1}(\theta)$
is an open set in ${\cal P}(\RR^{2d})$, and this concludes the
proof.

\Box

\begin{coro}\label{coro:useful2}

Let $(E,\Sigma)$ be a measurable space,  and $ \lambda \in E
\mapsto (\mu_{\lambda},\nu_{\lambda})\in ({\cal P}_2(\RR^{d}))^2$
a measurable function. We consider the function $\Psi$ defined by
(\ref{eq:psi}) and let $F$ be a closed set of $\RR^d$. Then, the
set
 $$\left\{(\lambda,x): (\{x\}\times F)\cap
 \Psi(\mu_{\lambda},\nu_{\lambda})\not= \emptyset\right\}$$ belongs to
$ \Sigma\otimes {\cal B}(\RR^d)$. In  particular,  if for all
$\lambda\in E$,
$\Pi^*_c(\mu_{\lambda},\nu_{\lambda})=\{\pi_{\lambda}\}$ is a
singleton, the set
$$\tilde{F}:=\left\{(\lambda,x): (\{x\}\times F)\cap
 supp(\pi_{\lambda})\not= \emptyset\right\}$$
is measurable.
\end{coro}

{\bf Proof.} Without loss of generality, we assume that $F$ is
nonempty. Let us first show that for any open set $\theta$ of
$\RR^{2d}$, the set
$$G=\{z\in \RR^{d}: (\{z\}\times F )\cap \theta \not =\emptyset\}$$
is open. Indeed, for $x\in G$ there exists  $y\in F$ and
$\varepsilon>0$ such that $B(x,\varepsilon)\times
B(y,\varepsilon)\subset \theta$. In particular, for all $z\in
B(x,\varepsilon)$ one has $(z,y)\in \theta$ and so
$B(x,\varepsilon) \subset G.$ By definition of measurability, the
set-valued mappings  $(\lambda,x)\to \{x\}\times F$ and
$(\lambda,x)\to \Psi(\mu_{\lambda},\nu_{\lambda})-(\{x\}\times F)$
are thus measurable.  The latter mapping being also closed valued,
we conclude that
$$\bigg\{(\lambda,x):
\left[\Psi(\mu_{\lambda},\nu_{\lambda})-(\{x\}\times F)\right]\cap
\{0\}\not = \emptyset\bigg\}$$ is a measurable set, which finishes
the proof.

\Box

\medskip

Let us now focus on the proof of Theorem \ref{theo:main}

\medskip

{\bf Proof of Theorem \ref{theo:main}} Since any $\sigma$-finite
measure is equivalent to a finite one, we can assume without loss
of generality that $m$ is finite.

For a fixed $k\geq 1$, we denote by $(A_{n,k})_{n\in\mathbb{Z}^d}$
the partition of $\RR^d$  in dyadic half-open rectangles of size
$2^{-dk}$, that is
$$A_{n,k}:=\prod_{i=1}^d \bigg[\frac{n_i}{2^{k}},\frac{n_i+1}{2^{k}}\bigg), \mbox{  where }n=(n_1,\dots,n_d)\in
\mathbb{Z}^d.$$
Consider  the sets $ B_{n,k}=\lbrace (\la,x)\in E\times \RR^d
:(\lbrace x\rbrace\times A_{n,k})\cap supp(\pi_\la)\neq
\emptyset\rbrace$. Notice that since $A_{n,k}=\bigcup_{j\in \NN}
\prod_{i=1}^d
\bigg[\frac{n_i}{2^{k}},\frac{n_i+1}{2^{k}}-\frac{1}{2^{k+j}}\bigg]$,
one has
$$  B_{n,k}= \bigcup_{j\in \NN} \left\{ (\la,x):\left(\lbrace
x\rbrace\times \prod_{i=1}^d
\bigg[\frac{n_i}{2^{k}},\frac{n_i+1}{2^{k}}-\frac{1}{2^{k+j}}\bigg]\right)\cap
supp(\pi_\la)\neq \emptyset\right\},$$ and so $B_{n,k}$ is
measurable thanks to Corollary \ref{coro:useful2}.

  Denote now by
$a_{n,k}\in A_{n,k}$ the ``center'' of the set, and
 define a $\Sigma\otimes {\cal B}(\RR^d)-$measurable function by
\begin{equation}
T^k(\la,x)=\sum_{n\in\mathbb{Z}^d}a_{n,k}\mathbf{1}_{B_{n,k}}(\la,x).
\end{equation}

For each $\la\in E$, let $\nu_\la^k$ be the discrete measure
defined by pushing forward  $\mu_\la$ through $T^k$, that is,
$$ \nu_\la^k(A)=\int\mathbf{1}_{T^k(\la,x)\in A}\ \mu_\la(dx),\ A\in\mathcal{B}(\RR^d).$$
 Denote also by  $\tilde{E}\in
\Sigma$ a measurable set with $m(\tilde{E}^c)=0$ and such that
 for all $\tilde{\lambda} \in \tilde{E}$,  $H(\mu_{\lambda},\nu_{\lambda},c)$ holds.

By hypothesis, for each $\lambda \in \tilde{E}$ we have that
\begin{equation}\label{1B=1A}
\mu_\la(dx) \mbox{ almost surely: }\mathbf{1}_{B_{n,k}}(\la,x)
=\mathbf{1}_{{\lbrace x:T_\la(x)\in A_{n,k}\rbrace}}.
\end{equation}
 where $T_{\lambda}$ has been defined in the statement of Theorem \ref{theo:main}.
This implies that
\begin{equation*}
\nu^k_\la(\{a_{n,k}\})=\int\mathbf{1}_{B_{n,k}}(\la,x)\mu_\la(dx)
=\mu_\la({\lbrace x:   T_\la(x) \in A_{n,k}\rbrace} )
=\nu_\la(A_{n,k})
\end{equation*}
by definition of $T_\la$.

We now check that $(T^k)_{k\in \NN}$ is a cauchy sequence in
$L^1(E\times \RR^d,m(d\la )\mu_{\la}(dx))$. Fix $k\leq k'$, and
for each $n\in \mathbb{Z}^d$ denote by $\{A_{n',k'}\}_{n'}$ the
unique partition of $A_{n,k}$ in dyadic rectangles of size $2^{-d
k'}$. We then have that
\begin{equation*}
\begin{split}
\int_E\int_{\RR^d} |T^k(\la,x)-T^{k'}& (\la,x)|  \mu_\la(dx)m(d\la) \\ & =\int_{E}\int_{\RR^d}\sum_{n\in\mathbb{Z}^d}\sum_{n':A_{n',k'}\subset A_{n,k}} \mathbf{1}_{B_{n',k'}}(\la,x) |a_{n,k}-a_{n',k'}|\mu_\la(dx)m(d\la)\\
&=\int_E\sum_{n\in\mathbb{Z}^d}\sum_{n':A_{n',k'}\subset A_{n,k}}|a_{n,k}-a_{n',k'}|\nu_\la(A_{n',k'}) m(d\la)\\
&\leq\int_E\sum_{n\in\mathbb{Z}^d}2^{-k}\sum_{n':A_{n',k'}\subset A_{n,k}}\nu_\la(A_{n',k'}) m(d\la)\\
&\leq\int_E\sum_{n\in\mathbb{Z}^d}
2^{-k}\nu_\la(A_{n,k}) m(d\la)\\
&\leq2^{-k}\int_E \nu_\la(\RR^d) m(d\la)= 2^{-k} m(E),
\end{split}
\end{equation*}
and the Cauchy property follows since $m(E)<\infty$.

\medskip

Let us denote by $T$ the limit in $L^1(E\times \RR^d,m(d\la
)\mu_{\la}(dx))$ of the sequence $T^k$. Theorem \ref{theo:main}
will be proved by verifying that for all $\la$ in a set of
$\Sigma$ of full $m$-measure set, one has $\pi_\la(dx,dy)=
\mu_\la(dx) \delta_{T(\la,x)}(dy)$. Hence, it is enough to check
that
$$\int\mathbb{I}_{C\times A_{n,k}}(x,T(\la,x))\mu_\la(dx)=\pi_\la\left( C\times A_{n,k}\right)$$
for any  semi-open rectangle $C$ with dyadic extremes and all
$n\in\mathbb{Z}^d,k\in \NN$. We have for $\la\in \tilde{E}$ and
any $j\in \NN$ that
\begin{equation}\label{piT}
\begin{split}
\Big| \pi_{\la} \left( C\times A_{n,k}\right) -
\int\mathbf{1}_{C\times A_{n,k}} & (x,T(\la,x)) \mu_\la(dx) \Big| \\
\leq  &  \left| \pi_{\la} \left(C\times A_{n,k}\right) -
\int\mathbf{1}_{C\times A_{n,k}}(x,T^j(\la,x))\mu_\la(dx)\right|
\\ & + \int\Big| \mathbf{1}_{C\times A_{n,k}}(x,T^j(\la,x))-
\mathbf{1}_{C\times A_{n,k}}(x,T(\la,x))\Big| \mu_\la(dx)\\
:= &  \Delta_j \ + \Delta'_j. \\
\end{split}
\end{equation}

We approximate $\mathbf{1}_{ A_{n,k}}$ by a
 Lipschitz continuous function $f_{\la,\varepsilon}$  such that \
 $\|f_{\la,\varepsilon}\|_{\infty}\leq 1$ and
$\mu_{\la}(\{y: f_{\la,\varepsilon}(y)\not = \mathbf{1}_{
A_{n,k}}\})\leq \varepsilon$ (this is possible thanks to
$H(\mu_\la,\nu_\la,c)$, a)). Hence,  the  second term $\Delta'_j$
on the r.h.s. of \eqref{piT} is bounded by
\begin{equation*}
4\varepsilon +L_{\la,\varepsilon} \int\left|
T^j(\la,x)-T(\la,x)\right|\mu_\la(dx),
\end{equation*}
where $L_{\la,\varepsilon}$ is the Lipschitz constant of
$f_{\la,\varepsilon}$. Since $\int\left|
T^j(\la,x)-T(\la,x)\right|\mu_\la(dx)$ converges in $L^1(m(d\la))$
to $0$,  there is a subsequence $T^{j_i}$ and a set $\hat{E}\in
\Sigma$ of full measure such that the convergence holds for all
$\la \in \hat{E}$. Consequently, for all $\lambda \in
\bar{E}:=\tilde{E}\cap \hat{E}$ we get that
$$\limsup_{i \to \infty} \Delta'_{j_i}= \limsup_{i \to \infty}\int\left| \mathbf{1}_{C\times A_{n,k}}(x,T^{j_i}(\la,x))-
\mathbf{1}_{C\times A_{n,k}}(x,T(\la,x))\right| \mu_\la(dx) \leq
2\varepsilon,$$ and since the l.h.s. does not depend on
$\varepsilon$, this means that $\lim_{i\to
\infty}\Delta'_{j_i}=0$.

 The proof will be achieved be verifying
that for fixed $\la\in \bar{E}$, one has $\Delta_{j}=0$ for all
large enough $j$. For such $\la$, fix a Borel set $D_\la$ of
$\RR^d$ of full $\mu_{\lambda}$ measure where \eqref{1B=1A} is
everywhere true. Then,
\begin{equation*}
\begin{split}
\int\mathbf{1}_{C\times A_{n,k}}(x,T^j(\la,x))\mu_\la(dx)= &
\int\mathbf{1}_{(D_{\la}\cap C)\times
A_{n,k}}(x,T^j(\la,x))\mu_\la(dx)\\ = &
\int\mathbf{1}_{(D_{\la}\cap C) \times A_{n,k}}\left(x,\sum_{m\in\mathbb{Z}^d}a_{m,j}\mathbf{1}_{B_{m,j}}(\la,x)\right)\mu_\la(dx)\\
=& \int\mathbf{1}_{(D_{\la}\cap C) \times A_{n,k}}\left(x,\sum_{m:a_{m,j} \in A_{n,k}}a_{m,j}\mathbf{1}_{A_{m,j}}(T_\la(x))\right)\mu_\la(dx).\\
\end{split}
\end{equation*}
Remark now that for all $j\geq k$, $y\in A_{n,k}
\Longleftrightarrow \sum_{m:a_{m,j} \in
A_{n,k}}a_{m,j}\mathbf{1}_{A_{m,j}}(y) \in A_{n,k}.$ Then, for all
$j\geq k$,
\begin{equation*}
\int\mathbf{1}_{C\times A_{n,k}}(x,T^j(\la,x))\mu_\la(dx)=
\int\mathbf{1}_{(D_{\la}\cap C)\times
A_{n,k}}(x,T_\la(x))\mu_\la(dx) =   \pi_{\la} \left( C\times
A_{n,k}\right).
\end{equation*}

\Box

\section{Application: strong coupling for orthogonal martingale measures}

We now develop an application of Theorem \eqref{theo:main}. Let
$(\Omega,{\cal F}, {\cal F}_t,\PP)$ be a filtered probability
space and consider $M$  an adapted orthogonal martingale measure
on $\RR_+\times\RR^d$ (in the sense of Walsh \cite{Walsh}). Assume
that its covariance measure has the form $q_t(da)dk_t$, where
$q_t(\omega,da)$ is a predictable random probability measure
 on $\RR^d$
 with finite second moment and $k_t$  a predictable increasing process.
Let us also consider another predictable random probability
measure $\hat{q}_t(\omega,da)$
 on $\RR^d$ with finite second moment.

 We want to construct in the same probability space a
second martingale measure with covariance measure
$\hat{q}_t(\omega,da)dk_t$, in such a way that in some sense, the
distance between the martingale measures is controlled by the
Wasserstein distance between their covariance measures. Recall
that this distance is defined for $\mu,\nu \in  {\cal P}_2(\RR^d)$
by
$$W_2^2(\mu,\nu)=\inf_{\pi<^{\mu}_{\nu}}I(\pi)$$
with  the quadratic cost $$c(x,y)=|x-y|^2.$$
 This distance makes the set
${\cal P}_2(\RR^d)$ a Polish space, and strengthens the weak
topology with the convergence of second moments (see
\cite{Rachev:98}).

\begin{theo} \label{theo:coupling}
In the previous setting,  assume moreover that
 $\PP(d\omega)dk_t(\omega)$ a.e.  $q_t$ has a density with respect to Lebesgue
 measure in $\RR^d$.

 Then,   there exists in $(\Omega,{\cal F}, {\cal F}_t,\PP)$  a martingale measure $\hat{M}$ on
$\RR_+\times\RR^d$ with covariance measure $\hat{q}_t(da)dk_t$,
such that for all $S>0$ and for every predictable function
$\phi:\Omega\times \RR_+\times \RR^d\to \RR$ that is Lipschitz
continuous in the last variable with $\EE\left(\int_0^S \int
\phi^2(s,a) \left(q_s(da)+\hat{q}_s(da)\right)
dk_s\right)<\infty$,  one has
\begin{equation}\label{eq:doobtransp}
\EE\left(\sup_{t\leq S}\left(\int_0^t \int \phi(s,a) M(ds,da)
-\int_0^t \int \phi(s,a) \hat{M}(ds,da) \right)^2 \right) \leq
\EE\left(\int_0^S L_s^2 \ W_2^2(q_s,\hat{q}_s)\ dk_s\right),
\end{equation}
where $L_s(\omega)$ is a measurable version of a Lipschitz
constant of $\phi(s,\omega,\cdot)$ and $W_2^2$ is the quadratic
Wasserstein distance in ${\cal P}_2(\RR^d)$.
\end{theo}

{\bf Proof.}  Since $q_t(\omega,da)$ has a density for almost
every $(t,\omega)$, assumption
 $\ H(q_t(\omega,da),\hat{q}_t(\omega,da),c)\ $ is satisfied. We can
 therefore
 apply Theorem \ref{theo:main} to $(E,\Sigma,m)=(\Omega\times
 \RR_+, {\cal P}red, \PP(d\omega)dk_t(\omega))$, where ${\cal
 P}red$ is the predictable $\sigma-$field with respect to ${\cal
 F}_t$. Then, there exists
 a  predictable mapping $T:\Omega\times\RR_+\times \RR^d:\to
 \RR^d$ that for $m$-almost every $(t,\omega)$
 pushes forward $q_t$ to $\hat{q}_t$. Moreover, for a.e. $(t,\omega)$, one has $$\int
 |a-T(t,\omega,a)|^2 q_t(\omega,da)=W_2^2(q_t,\hat{q}_t).$$

  On can thus define a martingale measure $\hat{M}$ by the stochastic integrals
 $$\int_0^t \int \psi(s,a) \hat{M}(ds,da):=\int_0^t \int \psi(s,T(s,a))
 M(ds,da)$$
for predictable simple functions $\psi$. Its  covariance measure
is by construction $\hat{q}_t(da) dk_t$, and by Doob's inequality,
the left hand side of \eqref{eq:doobtransp} is less than
\begin{align*}
\EE\Big( & \int_0^S \int |\phi(s,a)- \phi(s,T(s,a)) |^2
q_s(da)dk_s \Big) \leq \EE\left(\int_0^S L_s^2 \left(\int |a-
T(s,a) |^2 q_s(da)\right)dk_s \right) \\
 & \leq \EE\left(\int_0^S
L_s^2 \ W_2^2(q_s,\hat{q}_s) \ dk_s\right), \quad
 \mbox{ by definition of $T$ and of $W_2^2$. }
\end{align*}
\Box

\medskip

\begin{rem}
\begin{itemize}
\item[i)] The construction of strong couplings between orthogonal martingale
measures arises classically  in the literature, especially in
cases where  the martingale measure $M$ is a compensated Poisson
point measure  or a space-time white noise, for which the
covariance measures are deterministic  (cf. Grigelionis \cite{G},
El Karoui-Lepeltier \cite{EKL}, Tanaka \cite{T}, El
Karoui-M\'el\'eard \cite{EKM}, M\'el\'eard-Roelly \cite{MR},
Gu\'erin \cite{Gu}).
 A classical approach is to use the Skorokhod
representation theorem. This however prevents any hope to obtain
quantitative estimates related to the associated covariance
measures, what we have been able to do here thanks to the optimal
transport maps.
\item[ii)]  If the probability
space and the martingale measure $M$ are not fixed in advance, a
coupling  satisfying the estimate \eqref{eq:doobtransp} can be
constructed from and orthogonal martingale measure
$\tilde{M}(dt,da,da')$ on $\RR^+\times\RR^d\times\RR^d$ with
covariance measure $\pi_t(da,da')dk_t$, where $\pi_t$ is an
optimal transference plan between $q_t$ and $\hat q_t$.  Then,
$\tilde{M}(dt,da,\RR^d)$ and $\tilde{M}(dt,\RR^d,da')$ are indeed
two orthogonal martingale measures with the required covariances
and satisfying estimate \eqref{theo:coupling}. The question in
this situation is however how to construct such $\tilde{M}$.
\end{itemize}
\end{rem}

\end{document}